\RequirePackage{fix-cm}
\documentclass[smallextended]{svjour3}       % onecolumn (second format)
\smartqed  % flush right qed marks, e.g. at end of proof
\usepackage{graphicx}
\usepackage{amsmath, amssymb}
\begin{document}

\title{Regularized quasi-monotone method for stochastic optimization %\thanks{Grants or other notes
}


\author{V. Kungurtsev         \and
        V. Shikhman %etc.
}


\institute{Vaycheslav Kungurtsev \at
            Department of Computer Science, Faculty of Electrical Engineering \\ Czech Technical University in Prague \\ 121 35 Praha 2, Karlovo náměstí 13, Prague, Czech Republic \\ \email{kunguvya@fel.cvut.cz}         %  \\
           \and
           Vladimir Shikhman \at
              Department of Mathematics \\ Chemnitz University of Technology \\
Reichenhainer Str. 41, 
09126 Chemnitz, Germany \\
\email{vladimir.shikhman@mathematik.tu-chemnitz.de} 
}

\date{Received: date / Accepted: date}

\maketitle

\begin{abstract} We adapt the quasi-monotone method from~\cite{nesterov2015quasi} for composite convex minimization in the stochastic setting. For the proposed numerical scheme we derive the optimal convergence rate of  $\mbox{O}\left(\frac{1}{\sqrt{k+1}}\right)$ in terms of the last iterate, rather than on average as it is standard for subgradient methods. The theoretical guarantee for individual convergence of the regularized quasi-monotone method is confirmed by numerical experiments on $\ell_1$-regularized robust linear regression.

\end{abstract}

\section{Introduction}
\label{intro}
In the minimization of nonsmooth convex functions, typically, algorithms generate
a sequence of iterates using subgradients or estimates thereof. The convergence rates are then derived for some linear
combination of the iterates, rather than for the last estimate computed. 
Obtaining guarantees on the last iterate per se is often a challenging task. A significant contribution in that direction -- sometimes also refered to as individual convergence -- was given in~\cite{nesterov2015quasi}
with the quasi-monotone subgradient method. The corresponding analysis was simplified and extended
to solving minimization problems on decentralized networks in~\cite{liang2019distributed}.
In this paper we extend the work of~\cite{nesterov2015quasi} in two important directions,
first we consider a composite minimization problem with a simple additive function
(usually a regularizer), and second we consider the stochastic case. We develop the
Lyapunov-like analysis from~\cite{liang2019distributed} to handle the new elements and present
numerical experiments confirming the performance guarantees. We obtain the convergence rate of order $\mbox{O}\left(\frac{1}{\sqrt{k+1}}\right)$ in expectation, which is optimal for nonsmooth convex optimization. 

Let us briefly comment on the related literature. In~\cite{tao2020mirror} the authors introduce an adaptation of mirror
descent in order to attain the optimal individual convergence. They successively apply the latter for regularized nonsmooth learning problems in the stochastic setting. As shown in~\cite{tao2020extra}, the Nesterov's acceleration alternatively provides the individual convergence of projected subgradient methods as applied to nonsmooth convex optimization. Especially, the suggested methodology guarantees the regularization structure while keeping an optimal rate of convergence. Our contribution to individual convergence consists in theoretically justifying that also the initially proposed quasi-monotone subgradient method from~\cite{nesterov2015quasi} can be successively adjusted for composite  minimization in the stochastic setting.

\section{Regularized quasi-monotone method}

We consider the composite minimization problem
\begin{equation}
    \label{eq:opt}
    \min_{x} F(x)=\bar{f}(x)+g(x),
\end{equation}
where $\bar{f}, g: \mathbb{R}^n \rightarrow \mathbb{R} \cup \{+\infty\}$ are closed convex functions. Moreover, $$\bar{f}(x)=\mathbb{E}\left[f(x,\xi)\right]$$ for some $f$ closed
and convex in the first argument and $\xi$ is a sample from some random space $\Xi$. We assume that $\mbox{dom}(f(\cdot,\xi)) \subset \mbox{dom}(g)$ for a.e. $\xi$, and $\mbox{dom}(g)$ is closed. Usually, $\bar{f}$ plays the role of a loss function, whereas $g$ is used for regularization. In our setting, $f$ need not to be differentiable, but unbiased finite variance estimates of its subgradients, i.e. $w(x,\xi)\sim \nabla  f(x,\cdot)$ with $\mathbb{E}\left[w(x,\xi)\right] \in \partial \bar f(x)$, should be available. Here, we use $\nabla \bar f(x)$ to 
denote an element of the convex subdifferential $\partial \bar f(x)=\partial\mathbb{E}[f(x,\xi)]$, i.e.
\begin{equation}
    \label{eq:f-conv}
    \bar{f}(y) \geq \bar{f}(x) + \langle \nabla \bar f(x), y - x\rangle, \quad y \in \mbox{dom}(g).
\end{equation}
In addition, $g$ has to be simple. The latter means that we are able to find a closed-form solution for minimizing the sum of $g$ with some simple auxiliary functions. For that, we assume that for the effective domain of $g$ there exists a prox-function $\Psi:\mathbb{R}^n \rightarrow \mathbb{R} \cup \{+\infty\}$ w.r.t. an arbitrary but fixed norm $\|\cdot\|$. The prox-function $\Psi$ has to fulfil:
\begin{itemize}
    \item[(i)] $\Psi(x) \geq 0$ for all $x \in \mbox{dom}(g)$.%, and $\Psi(x_0)=0$ for some $x_0 \in \mbox{dom}(g)$. 
    \item[(ii)] $\Psi$ is strongly convex on $\mbox{dom}(g)$ with convexity parameter $\beta >0$, i.e. for all $x,y \in \mbox{dom}(g)$ and $\alpha \in [0,1]$ it holds:
    \[
     \Psi(\alpha x + (1-\alpha)y) \leq \alpha \Psi(x) + (1-\alpha) \Psi(y) - \frac{\beta}{2} \alpha (1-\alpha) \|x-y\|^2. 
    \]
    \item[(ii)] The auxiliary minimization problem
    \[
       \min_{x} \left\{\langle s,x\rangle + g(x) + \gamma \Psi(x)\right\} 
    \]
    is easily solvable for $s \in \mathbb{R}^n$ and $\gamma > 0$. 
    %Denote by $x_\gamma(s)$ its unique solution.
\end{itemize}
In our analysis, we consider that $g$ is strongly convex with convexity parameter $\sigma \geq 0$ w.r.t. the norm $\|\cdot\|$. Note that $\sigma=0$ corresponds to the mere convexity of $g$.

For stating our method, we choose a sequence  of positive parameters $(a_k)_{k \geq 0}$, which is used to average the subdifferential information of $f$. We set:
\[
  A_k = \sum_{\ell=0}^{k} a_\ell.
\]
Equivalently, it holds:
\begin{equation}
    \label{eq:Ak}
    A_{k+1}=A_k+a_{k+1}.
\end{equation}
Another sequence of positive parameters $(\gamma_k)_{k \geq 0}$ controls the impact of the prox-function $\Psi$. We assume:
\begin{equation}
    \label{eq:gamma-k}
 \gamma_{k+1} \geq \gamma_k, \quad k \geq 0.
\end{equation}
Now, we are ready to formulate the regularized quasi-monotone method for solving the composite minimization problem (\ref{eq:opt}):
\[
% \tag{DA}\label{eq:method}
   \begin{tabular}{|c|}
      \hline \\
      \begin{tabular}{c}
         {\bf Regularized Quasi-Monotone Method (RQM)} 
      \end{tabular}\\ \hline \\
      \begin{tabular}{ll}
            {\bf 0.} Initialize
            $\displaystyle x_{0} = \mbox{arg} \min_{x} \left\{ A_{0} g(x) + \gamma_{0}\Psi(x)\right\}$, $s_{-1}=0$.\\ \\
            {\bf 1.} Sample $\xi_k\sim \Xi$. \\ \\
            {\bf 2.} Compute $w\left(x_k,\xi_k\right)$ and set $\displaystyle s_{k}=s_{k-1} + a_{k} w\left(x_k,\xi_k\right)$. \\ \\
            {\bf 3.} Forecast 
            $\displaystyle x^+_{k} = \mbox{arg} \min_{x} \left\{ \left\langle s_{k}, x \right\rangle + A_{k+1} g(x) + \gamma_{k+1}\Psi(x)\right\}$.\\ \\
            {\bf 4.} Update 
            $\displaystyle x_{k+1} =  \frac{A_{k}}{A_{k+1}} x_k + \frac{a_{k+1}}{A_{k+1}} x^+_{k}$.\\
      \end{tabular}
      \\ \\ \hline
   \end{tabular}
\]
It is clear that iterates of (RQM) are convex combinations of forecasts:
\begin{equation}
  \label{eq:it-cv}
    x_{k} = \frac{1}{A_k} \left( a_0 x_0 + \sum_{\ell=1}^{k} a_{\ell} x^+_{\ell-1} \right).
\end{equation}

\section{Convergence analysis}

Before performing the convergence analysis of (RQM), let us deduce some useful properties of the following auxiliary function:
\begin{equation}
  \label{eq:phi-k}
    \varphi_{k}(s) = \max_{x\in\mathbb{R}^n} \left\{ \left\langle s, x \right\rangle - A_{k} g(x) - \gamma_{k}\Psi(x)\right\}, \quad s \in \mathbb{R}^n.
\end{equation}
Since $A_{k} g + \gamma_{k}\Psi$ is strongly convex with convexity parameter 
\begin{equation}
   \label{eq:mu-par}
   \mu_{k} = A_{k} \sigma + \gamma_{k} \beta, 
\end{equation}
the convex function $\varphi_{k}$ is differentiable and its gradient $\nabla \varphi_k$ is $\frac{1}{\mu_{k}}$-Lipschitz continuous. The latter property means:
\begin{equation}
    \label{eq:phi-lip}
    \varphi_{k}(s') \leq \varphi_{k}(s) + \langle \nabla \varphi_{k}(s), s'-s \rangle + \frac{1}{2 \mu_{k}} \|s'-s\|_*^2, \quad s,s' \in \mathbb{R}^n.
\end{equation}
Moreover, it holds:
\begin{equation}
    \label{eq:phi-plus}
  \nabla \varphi_{k}(-s_{k-1}) = x_{k-1}^+.
\end{equation}

Let us derive the convergence rate of (RQM). For that, we set:
\begin{equation}
  \label{eq:B-k}
    B_k = \frac{1}{2} \sum_{\ell=0}^{k} \frac{a_\ell^2}{\mu_{\ell}} \mathbb{E}\left\|w(x_\ell,\xi_{\ell})\right\|_*^2, \quad k\geq 0,
\end{equation}
where $\|\cdot\|_*$ denotes the dual norm of $\|\cdot\|$.
We shall denote, as standard, the filtration $\sigma$-algebra corresponding to the sequence of iterates as $\{\mathcal{F}_k\}$.
\begin{theorem} 
\label{thm:conv} Let $x_* \in \mbox{dom}(g)$ solve the composite optimization problem (\ref{eq:opt}), and the sequence $\left(x_k\right)_{k\geq 0}$ be generated by (RQM). Then, it holds for $k \geq 0$ that:
\begin{equation}
   \label{eq:cvg-f}
   \mathbb{E}\left[F(x_k)\right] - F(x_*) \leq \displaystyle
   \frac{\gamma_{k}}{A_k}\Psi(x_*) + \frac{B_k}{A_k}.
\end{equation}
%If additionally $g$ is nonnegative, we have:
%\begin{equation}
%   \label{eq:cvg-x}
%  \begin{array}{rcl}
%  \displaystyle
%  \frac{1}{2}\left\|x^+_{k-1} - x_*\right\|^2 &\leq& \displaystyle \frac{a_{k}}{\mu_{k}}\left(g(x_*)-\inf g(x)\right) + \frac{\gamma_{k}}{\mu_{k}}\Psi(x_*) + \frac{B_{k-1}}{\mu_{k}}.
%  \end{array}
%\end{equation}
\end{theorem}

\proof Let us define the stochastic Lyapunov function:
\[
   \begin{array}{rcl}
   V_k &=& A_k \left(F(x_k) - F(x_*)\right)  \displaystyle + \varphi_{k}(-s_k) + \langle s_k, x_*\rangle+ A_k g(x_*) - B_k.
   \end{array}
\]
We consider the expected difference:
\[
  \begin{array}{rcl}
      \mathbb{E}\left[V_{k+1}|\mathcal{F}_k\right] - V_{k} & = & \underbrace{A_{k+1}\left(\mathbb{E}\left[F(x_{k+1})|\mathcal{F}_k\right] - F(x_*)\right) - A_k \left(F(x_k) - F(x_*)\right)}_{\mbox{=I}} \\ \\
      &  & + \underbrace{\mathbb{E}\left[\varphi_{k+1}(-s_{k+1})|\mathcal{F}_k\right] - \varphi_{k}(-s_k)}_{\mbox{=II}} \\ \\
      & & + \underbrace{\mathbb{E}\left[\langle s_{k+1}, x_*\rangle |\mathcal{F}_k\right] + A_{k+1} g(x_*) - \langle s_k, x_*\rangle - 
A_k g(x_*)}_{\mbox{=III}} \\ \\
      & & \displaystyle \underbrace{-\mathbb{E}\left[B_{k+1}|\mathcal{F}_k\right]+ B_k}_{\mbox{=IV}}.
  \end{array}  
\]
Let us estimate the expressions I-IV from above.

{\bf Estimation of I:}
We split:
\[
   \begin{array}{rcl}
    \mbox{I} & = & \underbrace{A_{k+1}\left(\mathbb{E}\left[\mathbb{E}\left[f(x_{k+1},\xi)\right]|\mathcal{F}_k\right] - \mathbb{E}\left[f(x_*,\xi)\right]\right) - 
A_k \left(\mathbb{E}\left[f(x_k,\xi)\right] - \mathbb{E}\left[f(x_*,\xi)\right]\right)}_{=\mbox{I}_f} \\ \\
    && +\underbrace{A_{k+1}\left(\mathbb{E}\left[g(x_{k+1})|\mathcal{F}_k\right] - g(x_*)\right) - A_k \left(g(x_k) - g(x_*)\right)}_{=\mbox{I}_g}.
   \end{array}  
\]
Due to convexity of $f$, the definitions of $A_{k}$ and $x_{k}$, we obtain:
\[
    \begin{array}{rcl}
    \mbox{I}_f & \overset{(\ref{eq:Ak})}{=} & a_{k+1} \left(\mathbb{E}\left[\mathbb{E}\left[f(x_{k+1},\xi)\right]|\mathcal{F}_k\right] - \mathbb{E}\left[f(x_*,\xi)\right]\right) \\ \\ &&+ 
    A_k \left(\mathbb{E}\left[\mathbb{E}\left[f(x_{k+1},\xi)\right]|\mathcal{F}_k\right] - \mathbb{E}\left[f(x_k,\xi)\right]\right) \\ \\
    & \overset{(\ref{eq:f-conv})}{\leq} & a_{k+1} \mathbb{E}\left[\left\langle \nabla \bar f(x_{k+1}), x_{k+1} - x_*\right\rangle  |\mathcal{F}_k\right] + 
A_k \mathbb{E}\left[\left\langle \nabla \bar f(x_{k+1}), x_{k+1} - x_k\right\rangle |\mathcal{F}_k\right] \\ \\
    & \overset{\bf 4.}{=} & %\langle \nabla f(x_{k+1}), A_{k+1} x_{k+1} - A_k x_k - a_{k+1} x_*\rangle
     \mathbb{E}\left[\left\langle a_{k+1}\nabla  \bar f(x_{k+1}), x^+_{k} -  x_*\right\rangle|\mathcal{F}_k\right].
   \end{array}  
\]
By using convexity of $g$, it also follows:
\[
    \begin{array}{rcl}
    \mbox{I}_g & \overset{\bf 4.}{\leq} & \displaystyle A_{k+1}\left(\frac{A_k}{A_{k+1}}g(x_{k}) + \frac{a_{k+1}}{A_{k+1}}\mathbb{E}\left[g(x^+_{k})|\mathcal{F}_k\right] - g(x_*)\right) - A_k \left(g(x_k) - g(x_*)\right) \\ \\
    
    & \overset{(\ref{eq:Ak})}{=} & 
    a_{k+1}\left(\mathbb{E}\left[g(x^+_{k})|\mathcal{F}_k\right]-g(x_*)\right).
   \end{array}  
\]
Overall, we deduce:
\[
   \mbox{I} \leq  \mathbb{E}\left[\left\langle a_{k+1} \nabla \bar f(x_{k+1}), x^+_{k} -  x_*\right\rangle|\mathcal{F}_k\right]+ a_{k+1}\left(\mathbb{E}\left[g(x^+_{k})|\mathcal{F}_k\right]-g(x_*)\right).
\]

{\bf Estimation of II:}
First, in view of the definitions of $\varphi_k$, $A_k$, and $x^+_{k}$, we obtain:
\begin{equation}
 \label{eq:add1}
 \begin{array}{rcl}
    \mathbb{E}\left[\varphi_{k}(-s_k)| \mathcal{F}_k\right] & \overset{(\ref{eq:phi-k})}{\geq} & 
\mathbb{E}\left[\displaystyle \left\langle -s_k, x^+_{k} \right\rangle\right] - A_{k} \mathbb{E}\left[g(x^+_{k}) - \gamma_{k}\Psi(x^+_{k})|\mathcal{F}_k\right] \\ \\
    & \overset{(\ref{eq:Ak})}{=} & \mathbb{E}\left[\left\langle -s_k, x^+_{k} \right\rangle - A_{k+1} g(x^+_{k}) - \gamma_{k+1}\Psi(x^+_{k})\right] \\ \\ && +\mathbb{E}\left[ a_{k+1} g(x^+_{k}) + \left(\gamma_{k+1} - \gamma_{k}\right)\Psi(x^+_{k})|\mathcal{F}_k\right] \\ \\
    & \overset{\bf 3.}{=} & \mathbb{E}\left[\varphi_{k+1}(-s_k) + a_{k+1} g(x^+_{k}) + \left(\gamma_{k+1} - \gamma_{k}\right)\Psi(x^+_{k})|\mathcal{F}_k\right].
   \end{array}  
\end{equation}
Second, due to Lipschitz continuity of $\nabla \varphi_k$ and definitions of $s_k$ and $x_k^+$, we have:
\begin{equation}
 \label{eq:add2}
    \begin{array}{rcl}
    \mathbb{E}\left[\varphi_{k+1}(-s_{k+1})\right] & \overset{(\ref{eq:phi-lip})}{\leq} & \displaystyle \mathbb{E}\left[\varphi_{k+1}(-s_k) + \langle \nabla \varphi_{k+1}(-s_k), -s_{k+1}+s_{k} 
\rangle |\mathcal{F}_k\right] \\ \\ &&\displaystyle + \frac{1}{2 \mu_{k+1}} \mathbb{E}\left[\|-s_{k+1}+s_{k}\|_*^2|\mathcal{F}_k\right] \\ \\
     & \overset{{\bf 2.}, (\ref{eq:phi-plus})}{=} & \displaystyle \mathbb{E}\left[\varphi_{k+1}(-s_k) - \langle x^+_{k}, a_{k+1}w(x_{k+1},\xi_{k+1}) \rangle|\mathcal{F}_k\right] \\ \\ && \displaystyle + 
\frac{a_{k+1}^2}{2 \mu_{k+1}} \mathbb{E}\left[\|w(x_{k+1},\xi_{k+1}) \|_*^2|\mathcal{F}_k\right].
   \end{array}  
\end{equation}
By using these two auxiliary inequalities, we are ready to estimate:
\[
  \begin{array}{rcl}
    \mbox{II} & = & \mathbb{E}\left[\varphi_{k+1}(-s_{k+1}) - \varphi_{k}(-s_k)|\mathcal{F}_k\right] \\ \\
    & \overset{(\ref{eq:add1})}{\leq} & \mathbb{E}\left[\varphi_{k+1}(-s_{k+1}) - \varphi_{k+1}(-s_k) - a_{k+1} g(x^+_{k}) - \left(\gamma_{k+1} - \gamma_{k}\right)\Psi(x^+_{k})|\mathcal{F}_k\right] \\ \\
    & \overset{(\ref{eq:add2})}{\leq} & \displaystyle - \mathbb{E}\left[\langle a_{k+1} w(x_{k+1},\xi_{k+1}), x^+_{k} \rangle + \frac{a_{k+1}^2}{2 \mu_{k+1}} \|w(x_{k+1},\xi_{k+1}) \|_*^2|\mathcal{F}_k\right] \\ \\
    && - \mathbb{E}\left[a_{k+1} g(x^+_{k}) - \left(\gamma_{k+1} - \gamma_{k}\right)\Psi(x^+_{k})|\mathcal{F}_k\right].
   \end{array}  
\]

{\bf Estimation of III:}

The definitions of $s_k$ and $A_k$ provide:
\[
    \mbox{III} \overset{\bf 2.}{=}  \mathbb{E}\left[\langle a_{k+1} w(x_{k+1},\xi_{k+1}), x_*\rangle|\mathcal{F}_k\right] + a_{k+1}g(x_*).
\]

{\bf Estimation of IV:}

Here, we have:
\[
    \mbox{IV} \overset{(\ref{eq:B-k})}{=}  -\frac{a^2_{k+1}}{2\mu_{k+1}} \mathbb{E}\left[ \left\|w(x_{k+1},\xi_{k+1})\right\|_*^2|\mathcal{F}_k\right].
\]
Altogether, we can see that 
\[
\begin{array}{rcl}
  \mathbb{E}\left[V_{k+1}|\mathcal{F}_k\right] - V_k &\leq & 
\mathbb{E}\left[\langle a_{k+1} w(x_{k+1},\xi_{k+1}), x_*-x^+_k\rangle|\mathcal{F}_k\right] \\ \\ & & -\mathbb{E}\left[\left\langle a_{k+1} \nabla f(x_{k+1}), x^+_{k} -  x_*\right\rangle|\mathcal{F}_k\right] \\ \\ &&- \left(\gamma_{k+1} - \gamma_{k}\right)\Psi(x^+_{k}).
\end{array}
\]
Since $x^+_k$ is defined given $\mathcal{F}_k$, we have:
\[
\mathbb{E}\left[\langle a_{k+1} w(x_{k+1},\xi_{k+1}), x_*-x^+_k\rangle|\mathcal{F}_k\right] = \mathbb{E}\left[\langle a_{k+1} \nabla \bar f(x_{k+1}), x_*-x^+_k\rangle|\mathcal{F}_k\right].
\]
 By additionally using that the sequence $(\gamma_k)_{k \geq 0}$ is by assumption nondecreasing, and $\Psi(x) \geq 0$ for all $x \in \mbox{dom}(g)$, we obtain:
\[
   \mathbb{E}\left[V_{k+1}|\mathcal{F}_k\right] - V_k \leq 0.
\]
Hence, we get by induction and taking total expectations:
\begin{equation}
  \label{eq:add-v}  
  \mathbb{E}[V_k] \leq \mathbb{E}[V_0].
\end{equation}
It turns out that the expectation of $V_0$ is nonnegative. For that, we first estimate due to the choice of $x_0$:
\begin{equation}
   \label{eq:help1}
   \begin{array}{rcl}
   \varphi_0(-s_0) &\overset{(\ref{eq:phi-lip})}{\leq}& \displaystyle \varphi_{0}(0) + \langle \nabla \varphi_{0}(0), -s_0 \rangle + \frac{1}{2 \mu_{0}} \|s_0\|_*^2 \\ \\
   &\overset{\bf 0.}{=}& \displaystyle - a_0 g(x_0) - \gamma_0 \Psi(x_0) - \left\langle x_0, a_0 w(x_0,\xi_0) \right\rangle  \\ \\ &&+ \displaystyle \frac{a_0^2}{2\mu_{0}} \left\|w(x_0,\xi_0)\right\|_*^2.
   \end{array}
\end{equation}
This gives:
\begin{equation}
    \label{eq:v-0}   \begin{array}{rcl}
   \mathbb{E}[V_0] &=& \displaystyle A_0 \mathbb{E}\left[F(x_0) - F(x_*)\right]  + \mathbb{E}\left[\varphi_{0}(-s_0) + \langle s_0, x_*\rangle\right] \\ \\ &&+ A_0 g(x_*) - B_0 \\ \\ &\overset{(\ref{eq:f-conv})}{\leq}& \displaystyle a_0\left\langle \nabla  \bar f(x_0), x_0 \right\rangle  + \mathbb{E}[\varphi_{0}(-s_0)] + a_0 g(x_0) - B_0 \\ \\
   &\overset{(\ref{eq:help1}),(\ref{eq:B-k})}{\leq}& \displaystyle -
   \gamma_0 \Psi(x_0) \leq 0,
   \end{array}
\end{equation}
where again the last inequality is due to the assumptions on $\gamma_0$ and $\Psi$. Additionally, it holds by definition of $\varphi_k$:
\begin{equation}
  \label{eq:add-phi}
    \varphi_{k}(-s_k)  \geq \left\langle -s_k, x_* \right\rangle - A_{k} g(x_*) - \gamma_{k}\Psi(x_*).
\end{equation}
Hence, we obtain:
\[
   \begin{array}{rcl}
   A_k \mathbb{E}\left[\left(F(x_k) - F(x_*)\right)\right] &=& \displaystyle \mathbb{E}\left[V_k - \varphi_{k}(-s_k) - \langle s_k, x_*\rangle\right]- A_k g(x_*) + B_k \\ 
   &\overset{\scriptsize \begin{array}{c}
       (\ref{eq:add-v})   \\
       (\ref{eq:v-0}), (\ref{eq:add-phi})  
   \end{array}}{\leq} & \displaystyle
   \gamma_{k}\Psi(x_*) + B_k.
   \end{array}
\]
The assertion (\ref{eq:cvg-f}) then follows.

\qed

Now, let us show that the convergence rate of (RQM) derived in Theorem \ref{thm:conv} is optimal for nonsmooth optimization, i.e. it is of order $\mbox{O}\left(\frac{1}{\sqrt{k+1}}\right)$. For that, we exemplarily consider the following choice of control parameters:
\begin{equation}
    \label{eq:cpar1}
    a_k = 1, \quad \gamma_k = \sqrt{k+1}, \quad k \geq 0.
\end{equation}
We also assume that the subgradients' estimates of $f$ have uniformly bounded second moments, i.e. there exists $G >0$ such that
\begin{equation}
    \label{eq:lip-sg} 
    \mathbb{E}\left[\|w(x,\xi) \|_* \right] \leq G, \quad x \in \mbox{dom}(f).
\end{equation}  

\begin{corollary} 
\label{cor:conv1}
Let $x_* \in \mbox{dom}(g)$ solve the composite optimization problem (\ref{eq:opt}), and the sequence $\left(x_k\right)_{k\geq 0}$ be generated by (RQM) with control parameters from (\ref{eq:cpar1}). 
Then, it holds for all $k\geq 0$:
\begin{equation}
   \label{eq:cvg-f-c1}
   \mathbb{E}\left[F(x_k)\right] - F(x_*) \leq \displaystyle 
  \left(\Psi(x_*) + \frac{G^2}{\beta}\right) \frac{1}{\sqrt{k+1}}.
\end{equation}
%If additionally $g$ is bounded below, then all cluster points of $\left(x_k\right)_{k\geq 0}$ solve the composite optimization problem (\ref{eq:opt}).
\end{corollary}

\proof In order to obtain (\ref{eq:cvg-f-c1}), we estimate the terms in (\ref{eq:cvg-f}) which involve control parameters:
\[
    \frac{\gamma_k}{A_k} = \frac{\sqrt{k+1}}{k+1} = \frac{1}{\sqrt{k+1}},
\]
\[
 \begin{array}{rcl}
         \displaystyle  \frac{B_k}{A_k} &=& \displaystyle \frac{1}{2(k+1)}  \sum_{\ell=0}^{k} \frac{1}{\mu_{\ell}} \mathbb{E}\left[\left\|w(x_\ell,\xi_\ell)\right\|_*^2\right]
\overset{(\ref{eq:mu-par}), (\ref{eq:lip-sg})}{\leq}
\frac{G^2}{2 \beta (k+1)}  \sum_{\ell=0}^{k} \frac{1}{\sqrt{\ell+1}} \\ \\
 & \leq & \displaystyle \frac{G^2}{2 \beta (k+1)} \int_{-\frac{1}{2}}^{k + \frac{1}{2}} \frac{\mbox{d} \tau}{\sqrt{\tau+1}}  = \frac{G^2}{\beta (k+1)} \left( \sqrt{k+\frac{3}{2}} - \sqrt{\frac{1}{2}} \right) \\ \\ &\leq&
 \displaystyle \frac{G^2}{\beta \sqrt{k+1}}.
   \end{array}
\] \qed

%\begin{remark}
%Note that a high probability bound can be deduced by using the standard Hoeffding inequality. Namely, it holds:
%\[
%\begin{array}{l}
%\mathbb{P}_{\xi}\left[\sum\limits_{\ell=0}^k \left[f(x_\ell,\xi)+g(x_\ell)-f(x_*,\xi)+g(x_*)\right]  \ge 2\left(\Psi(x_*) + \frac{G^2}{\beta}\right) \sqrt{k+1}\right] \\ \qquad\qquad\le \exp\left\{-\frac{4(k+1)\left(\Psi(x_*) +\frac{G^2}{\beta}\right)^2}{2kF_M}\right\},
%\end{array}
%\]
%where $F_M = \sup\limits_{x,\xi} f(x,\xi)$, assuming such a bound %exists.\qed
%\end{remark}

We show that the convergence rate of (RQM) derived in Corollary \ref{cor:conv1} can be improved to $\mbox{O}\left(\frac{\ln k}{k}\right)$ if the regularizer $g$ turns out to be strongly convex. For that, consider the control parameters as follows:
\begin{equation}
    \label{eq:cpar2}
    a_k = 1, \quad \gamma_k = \ln (2 k +3), \quad k \geq 0.
\end{equation}

\begin{corollary} 
\label{cor:conv2}
Let $x_* \in \mbox{dom}(g)$ solve the composite optimization problem (\ref{eq:opt}), and the sequence $\left(x_k\right)_{k\geq 0}$ be generated by (RQM) with control parameters from (\ref{eq:cpar2}). Additionally, let $g$ be strongly convex with convexity parameter $\sigma >0$. Then, it holds for all $k\geq 0$:
\begin{equation}
   \label{eq:cvg-f-c2}
   \mathbb{E}\left[F(x_k)\right] - F(x_*) \leq \displaystyle 
  \left(\Psi(x_*) + \frac{G^2}{\sigma}\right) \frac{\ln (2 k +3)}{k+1}.
\end{equation}
%If additionally $g$ is bounded below, then all cluster points of $\left(x_k\right)_{k\geq 0}$ solve the composite optimization problem (\ref{eq:opt}).
\end{corollary}

\proof
In order to obtain (\ref{eq:cvg-f-c2}), we estimate the terms in (\ref{eq:cvg-f}) which involve control parameters:
\[
    \frac{\gamma_k}{A_k} = \frac{\ln (2 k +3)}{k+1},
\]
\[
 \begin{array}{rcl}
         \displaystyle  \frac{B_k}{A_k} &=& \displaystyle \frac{1}{2(k+1)}  \sum_{\ell=0}^{k} \frac{1}{\mu_{\ell}} \mathbb{E}\left[\left\|w(x_\ell,\xi_\ell)\right\|_*^2\right]
\overset{(\ref{eq:mu-par}), (\ref{eq:lip-sg})}{\leq}
\frac{G^2}{2 \sigma (k+1)}  \sum_{\ell=0}^{k} \frac{1}{\ell+1} \\ \\
 & \leq & \displaystyle \frac{G^2}{2 \sigma (k+1)} \int_{-\frac{1}{2}}^{k + \frac{1}{2}} \frac{\mbox{d} \tau}{\tau+1}  = \frac{G^2}{\sigma (k+1)} \left( \ln\left(k+\frac{3}{2}\right) - \ln{\frac{1}{2}} \right) \\ \\ &=&
 \displaystyle \frac{G^2}{\sigma}\cdot \frac{\ln(2k+3)}{k+1}.
   \end{array}
\] \qed

\section{Numerical experiments}
For our numerical illustration, let us consider linear regression with a robust Huber loss and $\ell_1$-regularization, i.e.
\begin{equation*}\label{eq:expobj}
\min\limits_{\mathbf{a},b} \,\,\sum\limits_{i=1}^N L_{\delta}\left(\mathbf{a}^T x_i+b-y_i\right)+\lambda\|(\mathbf{a},b)\|_1,
\end{equation*}
where
\begin{equation*}\label{eq:huber}
L_{\delta}(z)=\left\{\begin{array}{ll} 
\displaystyle \frac{1}{2} z^2 & \quad\text{for } |z|\le \delta,\\
\displaystyle \delta\left(|z|-\frac{1}{2}\delta\right) & \quad \text{otherwise}.
\end{array}\right.
\end{equation*}
Here, we expect the number $N$ of data samples to be large. The $\ell_1$-regularization
on the parameters encourages sparsity, i.e. most of the parameters to become zero. The Huber loss is a means of mitigating the impact of
outliers on the stability of the regression estimate, i.e. by enforcing linear as opposed to quadratic growth of the loss beyond the influence boundary $\delta$. We take the subgradients
\[
\partial L_{\delta}(z) \ni \left\{\begin{array}{ll} 
 z & \quad \text{for } |z|\le \delta, \\ \displaystyle
\delta \cdot \text{sign}(z) & \quad \text{otherwise}.
\end{array}\right.
\]
Denoting $x=(\mathbf{a},b)$ and choosing as prox-function $\Psi(x)=\frac{1}{2}\|x\|^2$, the subproblem in (RQM) admits an explicit solution:
\[
\begin{array}{rcl}
x^+ &=& \displaystyle \arg\min_x \left\{ \langle s,x\rangle + A \lambda \|x\|_1+\frac{\gamma}{2} \|x\|^2_2\right\} \\ \\&=& \displaystyle \text{sgn}\left(-\frac{s}{\gamma}\right)\max\left\{\left|-\frac{s}{\gamma}\right|- \frac{ A\lambda}{\gamma},0\right\}.
\end{array}
\]

\begin{center}
\begin{figure}[h]
\includegraphics[scale=0.5]{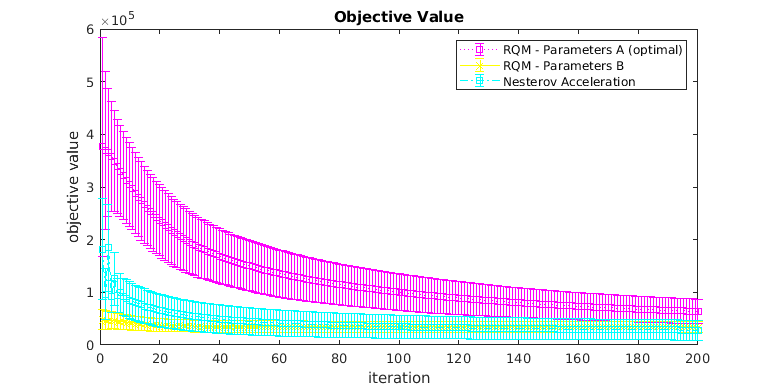}
\caption{Evolution of the objective error with the iterations.}
\label{fig:numerics}    
\end{figure}
\end{center}

To illustrate the performance of the algorithm, we first set up a synthetic data
profile. We let $\mathbf{a}\in\mathbb{R}^{10}$, $\delta=2$, and $N=10000$ and conduct the
following procedure:
\begin{enumerate}
\item Choose $4$ components of $\mathbf{a}$ to be nonzero. Randomly sample these components and $b$.
\item Choose $10000$ input samples $\{x_i\}$ uniformly in $[-5,5]$.
\item With probability $0.95$ generate $y_i \sim \mathcal{N}\left(\mathbf{a}^T x_i+b,1\right)$, and  otherwise $y_i \sim \mathcal{N}\left(\mathbf{a}^T x_i+b,5\right)$.
\item Run (RQM). 
\end{enumerate}
We ran one hundred trials of (RQM) in order to investigate the robustness and spread of the performance. We set $\lambda=0.1$. Note that the initial $x_0$ set by (RQM) is the
zero vector. First, we compare the parameter choice (\textbf{Parameters A}) as in (\ref{eq:cpar1}), i.e. $a_k=1$ and, thus, $A_k=k+1$, with $\gamma_k=\sqrt{k+1}$, to the choice of $a_k=k$ and, thus, $A_k = \frac{k(k+1)}{2}$, with the constant step-size $\gamma_k= 10$ (\textbf{Parameters B}). The trajectory of the objective value with the associated one standard deviation confidence interval is shown in Figure~\ref{fig:numerics}.

We also compare (RQM) to the stochastic regularized subgradient (SRSG) with Nesterov's extrapolation from~\cite{tao2020extra}. Their, by choosing control parameters $$\theta_k = \frac{2}{k+1}, \quad \gamma_k = (k+1)^{3/2},$$ the authors iterate:
\[
\begin{array}{l}
 y_{k} = \hat x_{k}+\theta_k\left(\theta^{-1}_{k-1}-1\right)\left(\hat x_k- \hat x_{k-1}\right), \\ \\ \displaystyle
\hat x_{k+1} = \mbox{arg} \min_{x} \left\{ \left\langle w(y_k,\xi_k), x \right\rangle +  g(x) + \gamma_{k}\Psi(x- y_k)\right\}.
\end{array}
\]
The explicit solution of the latter subproblem in our context is
\[
\begin{array}{rcl}
 \hat x &=& \displaystyle \arg\min_x \left\{ \langle w,x\rangle + \lambda \|x\|_1+\frac{\gamma}{2} \|x - y\|^2_2\right\} \\ \\&=& \displaystyle \text{sgn}\left(y-\frac{w}{\gamma}\right)\max\left\{\left|y-\frac{w}{\gamma}\right|- \frac{\lambda}{\gamma},0\right\}.
\end{array}
\]
Note that in Figure \ref{fig:numerics} we report the objective value on $\hat x_k$, as this is what the theoretical convergence guarantees in~\cite{tao2020extra} are derived for.

Overall, we conclude that while all the methods appear to be convergent, clearly the parameter choices suggested by the theory are a worst-case bound, and more aggressive parameter choices appear to work fine in practice.

\bibliographystyle{spmpsci}      % mathematics and physical sciences
\bibliography{refs}   % name your BibTeX data base

\end{document}